\documentclass[12pt]{article}
\input{xypic.tex} \xyoption{all}
\usepackage{amsthm,amssymb,hyperref}
\title{Spectral curves and Nahm transform for doubly-periodic instantons}
\author{Marcos Jardim \\ University of Pennsylvania \\ Department of Mathematics \\
Philadelphia, PA 19104-6395 USA \\ jardim@math.upenn.edu}

\newcommand{\seta}{\rightarrow} \newcommand{\oo}{{\cal O}}
\newcommand{\torus}{T\times\cpx} \newcommand{\tproj}{T\times\proj}
\newcommand{\dual}{\hat{T}} \newcommand{\dproj}{\dual\times\proj}
\newcommand{\del}{\overline{\partial}} \newcommand{\delad}{\overline{\partial}^*}
\newcommand{\dtorus}{\hat{T}\times\cpx} \newcommand{\real}{ \mathbb{R}} 
\newcommand{\cpx}{\mathbb{C}} \newcommand{\proj}{\mathbb{P}^1}
\newcommand{\call}{{\cal L}} \newcommand{\ee}{{\cal E}} \newcommand{\as}{\pm\xi_0} 

\newtheorem{thm}{Theorem} 
\newtheorem{prop}[thm]{Proposition} 

\begin{document}
\maketitle

\begin{abstract}
We explore the role played by the spectral curves associated with
Higgs pairs in the context of the Nahm transform of
doubly-periodic instantons previously defined by the author. More
precisely, we show how to construct a triple consisting of an
algebraic curve plus a line bundle with connection over it from 
a doubly-periodic instanton, and that these coincide with the 
Hitchin's spectral data associated with the Nahm transformed 
Higgs bundle.
\end{abstract}

\baselineskip18pt \newpage

\section{Introduction} \label{intro}

In \cite{J1,J2}, we have established a correspondence between
instantons on $\real^4$ which are periodic in two directions
(so-called {\em doubly-periodic instantons}) and certain singular
Higgs pairs over a 2-dimensional torus. On the other hand, Hitchin
has shown that Higgs pairs are equivalent to a pair consisting of
an algebraic curve (the {\em spectral curve}) in the total space
of the cotangent bundle plus a ``line bundle'' over it.

In this third installment of the series initiated by \cite{J1,J2},
we shall explore the relation between Hitchin's spectral
data and the Nahm transform for doubly-periodic instantons 
defined in the previous papers.

The motivation comes from the Hitchin's work on monopoles
\cite{H1}. He has shown that monopoles on $\real^3$ (that is,
instantons on $\real^4$  which are translation invariant in one
direction) are equivalent to certain singular solutions of Nahm's
equations and to pair consisting of an algebraic curve plus a line
bundle over it.

The paper is organized as follows. In section \ref{nt} we briefly
review the Nahm transform of doubly-periodic instantons discussed
in the previous papers. We then show how to construct the spectral
data from the instanton (section \ref{instspec}) and from the
Higgs bundle (section \ref{hitspec}). The main purpose of this paper 
is to prove that these two sets of spectral data coincide when the 
Higgs bundle is the Nahm transform of a doubly-periodic instanton; 
this is done in section
\ref{match}.

\paragraph{Acknowledgements.}
This work has been partially funded by the Brazilian Ministry of
Science and Technology via CNPq. I am grateful to Nigel Hitchin
for his help on the early stages of this project.


\section{Nahm transform} \label{nt}

Let $T$ be a complex torus, and consider the product $\torus$
equipped with the product flat metric. Let $E\seta\torus$ be a
rank 2 complex vector bundle with an irreducible $SU(2)$ instanton
(i.e. anti-self-dual) connection $A$ such that $|F_A|=O(|w|^{-2})$, 
where $w$ is a complex coordinate on $\cpx$. As usual, its total 
energy is given by:
$$ \int_{\torus} |F_A|^2 = 8 \pi^2 k $$
where $k$ is a positive integer, the instanton number.

As it was shown in \cite{BJ}, the toroidal components of the
connection $A$ have a well-defined limit as $r\seta\infty$ given
by a constant flat connection $\Gamma$ over $T$. General theory
tells us that a constant flat connection on a bundle $S\seta T$
determines uniquely a holomorphic structure on $S$. Moreover, $S$
must split, holomorphically, as the sum of two flat line bundles,
i.e. $S=\xi_0 \oplus -\xi_0$, uniquely up to $\pm1$. Here,
$\as$ are seen as points in $\dual$, the torus dual to $T$.

The following result has been proved in \cite{J1,J2}:

\begin{thm} \label{nahmthm}
The {\em Nahm transform} is a bijective correspondence between the
following objects:
\begin{itemize}
\item irreducible $SU(2)$ instanton connections on $E\seta\torus$ with
fixed instanton number $k$ and asymptotic state $\xi_0$; and
\item admissible $U(k)$ solutions of the Hitchin's equations over $\dual$,
such that the Higgs field has at most simple poles at
$\as\in\dual$, with semisimple residues of rank $\leq2$ if $\xi_0$
is an element of order 2 in the Jacobian of $T$, and rank $\leq1$
otherwise.
\end{itemize} \end{thm}

For the purpose of this paper, it is enough to recall one way of
the above correspondence. For simplicity, let us assume that
$\xi_0\neq-\xi_0$.

Recall that $\dual$ parametrises the set of line bundles with flat
connection on $T$. Indeed, given a point $\xi\in\dual$, we denote
by $L_\xi$ the trivial line bundle with the constant connection
$\omega_\xi=i\xi$. Let $\pi_1:\torus\seta T$ be the obvious
projection, and consider the twisted bundles $\ee(\xi) =
\ee\otimes\pi_1^*L_\xi$ with the corresponding instanton
connection:
$$ A_{\xi}=A\otimes{\rm Id} + {\rm Id}\otimes\omega_{\xi} $$
For each $\xi\neq\as$ we have the elliptic complex
\begin{equation} \label{monad}
0 \seta L^2_2(\Omega^0E(\xi))
\stackrel{\del_{A_\xi}}{\longrightarrow}
   L^2_1(\Omega^{0,1}E(\xi)) \stackrel{-\del_{A_\xi}}{\longrightarrow}
   L^2(\Omega^{0,2}E(\xi)) \seta 0
\end{equation}
whose 0th and 2nd cohomologies vanish, while $H^1(\xi)$ (which
coincides with the cokernel of the coupled Dirac operator
$D_{A_\xi}$
\footnote{Recall that on a K\"ahler manifold the Dirac operator
$D$ can be written as $D=\del-\delad$.}
) has complex dimension $k$. Thus, the cohomology of
the above monad defines a rank $k$ vector bundle
$V\seta(\dual\setminus\{\as\})$, with fibres
$V_\xi=H^1(\xi)={\rm ker}D_{A_\xi}^*$, plus an unitary
connection $B$ obtained as follows. Let $\cal H$ be the trivial
Hilbert bundle over $\dual\setminus\{\as\}$ with fibres given 
by $L^2_1(\Omega^{0,1}E(\xi))$. Then $V$ is a naturally a
subbundle of $\cal H$, and we denote by $P$ the projection 
${\cal H} \seta V$ and by $\iota_{V\seta{\cal H}}$ the inclusion 
$V\hookrightarrow {\cal H}$.
We define:
$$ \nabla_B = P \circ d \circ \iota_{V\seta{\cal H}} $$
where $d$ denotes the trivial connection on $\cal H$.

A Higgs field $\Phi\in{\rm End}(V)\otimes K_{\dual}$ is defined as follows.
Let $\psi$ be a section of $V$, i.e. for each $\xi\in\dual\setminus\{\as\}$, 
$\psi[\xi]\in{\rm ker}D^*_{A\xi}$. For a fixed $\xi$, the Higgs 
field will act on $\psi[\xi]$ by multiplication by the plane 
coordinate $w$ composed with projection to ${\rm ker}D^*_{A_\xi}$:
\begin{equation} \label{higgs.gt}
(\Phi(\psi))[\xi] \ = \ \frac{1}{\sqrt{2}}P_{\xi}(w\psi[\xi])d\xi
\end{equation}
where $P_{\xi}:L^2_1(\Omega^{0,1}E(\xi)) \seta {\rm ker}D_{A_\xi}^*$  
is the natural projection operator.

As it was shown in \cite{J2}, the pair $(B,\Phi)$ satisfy Hitchin's
equations \cite{H2}.

\paragraph{Kobayashi-Hitchin correspondence for doubly-periodic
instantons.}
It is also useful to recall that given an doubly-periodic instanton
connection $A$ on the bundle $E\seta\torus$ there is a holomorphic
bundle $\ee\seta\tproj$, so callec {\em instanton bundle} such that:
\begin{itemize}
\item $\ee|_{T\times(\proj\setminus\infty)} = (E,\del_A)$;
\item $c_1(\ee)=0$ and $c_2(\ee) = k$;
\item $\ee|_{T\times\infty} = L_{\xi_0}\oplus L_{-\xi_0}$.
\end{itemize}
Instanton bundles also satisfy a certain stability condition, 
but that is not relevant here; see \cite{BJ} for the complete 
statement.

It is then easy to show that $h^0(\ee(\xi))=h^2(\ee(\xi))=0$, while 
$H^1(\ee(\xi))$ can be identified with ${\rm ker}D_{A_\xi}^*$ \cite{J2}.

Furthermore, every instanton bundle is {\em generically fibrewise semistable}, 
that is $\ee|_{T_p}$ is semistable for generic $p\in\proj$. The instanton
bundle is said to be {\em fibrewise semistable} if $\ee|_{T_p}$ is semistable 
for every $p\in\proj$. It is {\em regular} if it is fibrewise semistable
and $h^0(\rm{End}(\ee|_{T_p}))=2$ for all $p\in\proj$
\footnote{Recall that every semistable, rank 2 vector bundle over an eliptic
with trivial determinant either splits as a sum of flat line bundles or it is
the unique nontrivial extension of a flat line bundle of order 2 by itself. Such
bundle is regular if it is not the sum of flat lines bundle of order 2.}. 

The instanton $A$ is said to be regular if the corresponding 
instanton bundle is regular. As we will see below, regular instantons form
a Zariski open subset of the moduli space of  doubly-periodic instanton.

Finally, it is important to remind how the Higgs field $\Phi$
can be constructed out of this holomorphic data. 

Start by fixing two sections $s_0$ and $s_\infty$ generating 
$H^0(\proj,\oo_{\proj}(1))$, such that $s_0$ vanishes at $0\in\cpx$ and 
$s_\infty$ vanishes at the point added at infinity. For each 
$\xi\in\dual\setminus\{\as\}$, 
we define the map:
\begin{equation} \begin{array}{rcl}
H^1(\tproj,\ee(\xi))\times H^1(\tproj,\ee(\xi)) &
\stackrel{\Psi_\xi}{\longrightarrow} & H^1(\tproj,\tilde{\ee}(\xi)) \\
(\alpha,\beta) & \mapsto & \alpha\otimes s_0-\beta\otimes s_\infty
\end{array} \end{equation}
If $(\alpha,\beta)\in{\rm ker}\Psi_\xi$, we define 
$\Phi\in{\rm End}(H^1(\tproj,\ee(\xi)))={\rm End}(V_\xi)$ as
follows:
\begin{equation} \label{alt.hig2}
\varphi_\xi(\alpha)=\beta
\end{equation}
To check that the two definitions coincide at $\cpx=\proj\setminus\infty$,
just note that $\frac{s_0(w)}{s_\infty(w)}=w$ for any trivialisation of 
$\oo_{\proj}(1)$ with local coordinate $w$.

With this formulation, it is not difficult to show that the eigenvalues 
of the Higgs field $\Phi$ have at most simple poles at $\as$. Moreover, 
the residues of $\Phi$ are semisimple and have rank $\leq2$ if $\xi_0$ 
is an element of order 2 in the Jacobian of $T$, and rank $\leq1$ otherwise;
see \cite{J2} for the details.


\section{The instanton spectral data} \label{instspec}

Our first step is to construct a complex curve
$S\hookrightarrow\dtorus$ associated to a doubly-periodic
instanton $A$. 

Let $D_{A_\xi(w)}^*$ denote the restriction of the coupled
Dirac operator $D_{A_\xi}$ to the torus $T_w$. We define:
\begin{equation} \label{def.spec}
S = \{ (\xi,w) \in \dtorus \ | \ {\rm ker} \{ D_{A_\xi(w)}^* \} \neq 0 \}
\end{equation}
Since $D_{A_\xi(w)} = \del_{A_\xi}|_{T_w} - \del^*_{A_\xi}|_{T_w}$, it is easy to see
that:
$$ {\rm ker} \{ D_{A_\xi(w)}^* \} = H^1(T_w,E(\xi)|_{T_w}) = H^1(T_w,\ee(\xi)|_{T_w}) $$
Note also that $S$ can be compactified to a curve $\overline{S}\hookrightarrow\dproj$ 
by adding the two points $(\as,\infty)$ corresponding to the asymptotic states.

Assuming that the instanton bundle is fibrewise semistable, we conclude that  
$\overline{S}$ is a branched double cover of $\proj$; the branch points correspond
to those $w\in\cpx$ such that $\ee(\xi)|_{T_w}$ is an extension of the trivial line 
bundle by itself. 

On the other hand, index theorem tells us that $\overline{S}$ is a $k$-fold cover of $\dual$. 
Hence there are $4k$ branch points, and the genus of $\overline{S}$ is $2k-1$. Moreover, 
all spectral curves belong to the linear system  $| k\cdot[\dual] + 2\cdot[\proj] | \subset \dproj$. 
The curve $\overline{S}$ is smooth provided $A$ is regular.

\paragraph{Line bundle with connection.}
Let $\pi_1:\dproj\seta\dual$ and $\pi_2:\dproj\seta\proj$ be the natural 
projection maps; we will also use $\pi_1$ and $\pi_2$ to denote the projections
$\overline{S}\seta\dual$ and $\overline{S}\seta\proj$.  
 
To each $s\in\overline{S}$, we attach the vector space:
\begin{equation} \label{spec.bdl1}
{\cal L}_s = {\rm ker} \left\{ D_{A_{\pi_1(s)}(\pi_2(s))}^* \right\} 
           = H^1(T_{\pi_2(s)},\ee(\pi_1(s))|_{T_{\pi_2(s)}}) 
\end{equation}
If $\ee$ is only fibrewise semistable, then ${\cal L}$ is only a coherent sheaf
on the (singular) spectral curve. However, when the instanton bundle is regular 
${\cal L}$ becomes a line bundle.

So now let us assume that $A$ is a regular doubly-periodic instanton, and consider the bundle
$\pi_1^*{\cal H} \seta S$. There is a bundle map $T : \pi_1^*{\cal H} \seta {\cal L}$, 
which is given by the following composition on each fibre:
\begin{equation} \label{compo}
L^2_1(\Omega^{0,1}E(\pi_1(s))) \stackrel{P}{\seta} {\rm ker} \left\{ D_{A_{\pi_1(s)}}^* \right\}  
\stackrel{r}{\seta} {\rm ker} \left\{ D_{A_{ \pi_1(s) }(\pi_2(s))}^* \right\}
\end{equation}
where $r$ denotes the restriction map. Let $\iota_{{\cal L}\seta{\cal H}}$ denotes 
the inclusion \linebreak ${\cal L} \hookrightarrow \pi_1^*{\cal H}$, which makes sense in 
terms of distributions.  A connection $\Gamma$ on the line bundle ${\cal L} \seta S$ 
is defined by:
\begin{equation} \label{conn1}
\nabla_\Gamma = T \circ \pi_1^*d \circ \iota_{{\cal L}\seta{\cal H}} 
\end{equation}


\section{Hitchin's spectral data} \label{hitspec}

We now look at the other side of the correspondence in theorem
\ref{nahmthm} and review Hitchin's construction of spectral curves
associated to Higgs bundles \cite{H3}.

Recall that $V\seta\dual\setminus\{\as\}$ is a rank $k$ vector bundle,
and $\Phi$ is an endomorphism valued $(1,0)$-form with simple
poles at $\as$. So, for any fixed $\xi\in\dual\setminus\{\as\}$,
$\Phi[\xi]$ is a $k\times k$ matrix and one can compute its $k$
eigenvalues. As we vary $\xi$, we get  a $k$-fold covering,
possibly branched, of $\dual\setminus\{\as\}$ inside
$\dtorus$. This {\em curve of eigenvalues} is what we want
to define as our {\em Higgs spectral curve}; more precisely:
\begin{equation} \label{h.spec}
C=\left\{ (\xi,w)\in\dtorus\ |\ {\rm det}(\Phi[\xi]-w\cdot{\rm I}_k)=0 \right\}
\end{equation}
In other words, $C$ is the set
of points $(\xi,w)\in\dtorus$ such that $w$ is an eigenvalue of
the endomorphism $\Phi[\xi]:V_\xi\seta V_\xi$. 

Since we are assuming that $\Phi$ has simple poles at $\as$, 
the curve $C\hookrightarrow\dtorus$ can be compactified to a curve 
$\overline{C}\hookrightarrow\dproj$ by adding the 
points $(\as,\infty)$.  

The following proposition is a familiar fact from the theory of 
Higgs bundles. 

\begin{prop} \label{smooth}
If $\xi_0\neq-\xi_0$, the spectral curve associated to a generic
Higgs bundle $(V,B,\Phi)$ is smooth.
\end{prop}

Note that if $\xi_0=-\xi_0$, then all spectral curves have 
a double-point at $(\as,\infty)$, but are generically smooth 
elsewhere.

\paragraph{Defining the spectral bundle.}
As before, we will denote the projections $\overline{C}\seta\dual$ and
$\overline{C}\seta\proj$ by $\pi_1$ and $\pi_2$. We define a coherent 
sheaf ${\cal N}$on $\overline{C}$ with stalks given by:
\begin{equation} \label{spec.bdl2}
{\cal N}_c = {\rm coker}\left\{ \Phi[\pi_1(c)] - \pi_2(c) \cdot {\rm Id}_k \right\}
\end{equation}
i.e. the dual of the $\pi_2(c)$-eigenspace of $\Phi[\pi_1(c)]$
Generically, one expects the eigenvalues to be distinct, so that
${\cal N}$ becomes a line bundle over the smooth curve $\overline{C}$.

Assuming that Higgs bundle $(V,B,\Phi)$ is generic, we define a connection
$\Lambda$ on the line bundle ${\cal N} \seta C$. First note that
$\cal N$ is naturally a subbundle of $\pi_1^*V$; let $\iota_{{\cal N}\seta V})$ 
be the inclusion and $E : \pi_1^*V \seta {\cal N}$ the fibrewise projection. 
We define: 
\begin{equation} \label{conn2}
\nabla_\Lambda = E \circ \pi^*_1\nabla_B \circ \iota_{{\cal N}\seta V}
\end{equation}


\section{Matching the spectral data} \label{match}

We are finally in a position to state and prove the main result of
this paper:

\begin{thm} \label{specthm}
If $(V,B,\Phi)$ is the Nahm transform of a regular instanton $(E,A)$, 
then the instanton spectral data $(\overline{S},\call,\Gamma)$ is 
equivalent to the Higgs spectral data $(\overline{C},{\cal N},\Lambda)$, 
in the sense that the curves $S$ and $C$ coincide pointwise and 
there is a natural isomorphism $\call\seta{\cal N}$ preserving 
the connections.
\end{thm}

\begin{proof}
Clearly, both spectral curves already have the points
$(\as,\infty)$ in common. So let $\xi\neq\as$ and suppose that
$\alpha$ is an eigenvector of $\Phi[\xi]$  with eigenvalue
$\epsilon<\infty$. In particular, the point
$(\xi,\epsilon)\in\dtorus$ belongs to the Higgs spectral curve
$C$. By definition, we have:
$$ \Phi[\xi](\alpha)=\epsilon\cdot\alpha \ \ \ \Rightarrow \ \ \
   \alpha\otimes(s_0-\epsilon\cdot s_\infty)=0 $$

Clearly, $s_\epsilon=s_0-\epsilon\cdot s_\infty$ is a holomorphic
section in $H^0(\proj,\oo_{\proj}(1))$ vanishing at
$\epsilon\in\proj\setminus\{\infty\}$. Therefore it induces the
following exact sequence:
$$ 0 \seta \ee(\xi) \stackrel{\otimes s_\epsilon}{\seta}
   \widetilde{\ee}(\xi) \seta \widetilde{\ee}(\xi)|_{T_\epsilon} \seta 0 $$
which in turn induces the cohomology sequence:
\begin{equation} \label{sqc.xx} \begin{array}{cccccc}
0 & \seta & H^0(T_\epsilon,\widetilde{\ee}(\xi)|_{T_\epsilon}) & \seta & & \\
  & \seta & H^1(\tproj,\ee(\xi)) & \stackrel{\otimes s_\epsilon}{\seta} &
            H^1(\tproj,\widetilde{\ee}(\xi)) & \stackrel{r}{\seta} \\
  & \stackrel{r}{\seta} & H^1(T_\epsilon,\widetilde{\ee}(\xi)|_{T_\epsilon}) &
\seta & 0 &
\end{array} \end{equation}
Thus $\alpha\in{\rm ker}(\otimes s_\epsilon)=
H^0(T_\epsilon,\widetilde{\ee}(\xi)|_{T_\epsilon})=
H^0(T_\epsilon,\ee(\xi)|_{T_\epsilon})=H^1(T_\epsilon,\ee(\xi)|_{T_\epsilon})^*$.

In particular, $H^1(T_\epsilon,\ee(\xi)|_{T_\epsilon})={\rm ker}\{D_{A_\xi(w)}^*\}$ is
non-empty, hence $(\xi,\epsilon)\in\dtorus$ also belongs to the instanton spectral curve 
$S$. The same argument clearly provides the converse statement. Thus the curves $C$ and 
$S$ must coincide pointwise.

It also follows from the cohomology sequence (\ref{sqc.xx}) that
the dual of the $\epsilon$-eigenspace of $\Phi[\xi]$ is exactly
$H^1(T_\epsilon,\widetilde{\ee}(\xi)|_{T_\epsilon})=
H^1(T_\epsilon,\ee(\xi)|_{T_\epsilon})$. In other words, there
are canonical identifications between the fibres 
${\cal N}_{(\xi,\epsilon)}$ and $\call_{(\xi,\epsilon)}$, and the line bundles 
are isomorphic.

Finally, let us check that the connection $\nabla_\Gamma$ and 
$\nabla_\Gamma$ also coincide. Noting that the projection
$E:\pi_1^*V\seta {\cal N}={\cal L}$ is just the restriction map 
\linebreak $r:{\rm ker} \left\{ D_{A_{\pi_1(s)}}^* \right\}  
\seta {\rm ker} \left\{ D_{A_{ \pi_1(s) }(\pi_2(s))}^* \right\}$
on each $s\in S=C$, it is easy to see that $T=E \circ \pi_1^*P$
Therefore, we have:
\begin{eqnarray*}
\nabla_\Gamma & = & T \circ \pi_1^*d \circ \iota_{{\cal N}\seta{\cal H}} = \\
& = & E \circ \left( \pi_1^*P \circ \pi_1^*d \circ \iota_{V\seta{\cal H}} \right) \circ \iota_{{\cal N}\seta V} =
   E \circ \pi_1^*\nabla_B \circ \iota_{{\cal N}\seta V} = \nabla_\Lambda
\end{eqnarray*} 
\end{proof}

\vskip12pt

\noindent {\bf Remark 1:} 
More generally, the above argument shows that the pairs $(\overline{S},\call)$ 
and $(\overline{C},{\cal N})$ also coincide when $A$ is fibrewise semistable.

\vskip12pt

\noindent {\bf Remark 2:}
Cherkis and Kapustin used a similar argument to establish the
analogous result for periodic monopoles \cite{CK}. More precisely,
they considered monopoles on $S^1\times\real^2$, so that the Nahm 
transformed object is a Higgs pair on $S^1\times\real$. Each of these
objects can be associated to a spectral pair consisting of an
algebraic curve on $\real^2\times(\real^2\setminus\{0\})$ plus a
line bundle over it. If the Higgs pair is the Nahm transform of a
periodic monopole, Cherkis and Kapustin have shown that both spectral
data coincide.


\section{Relation with Fourier-Mukai transform}

The instanton spectral pair $(\overline{S},\call)$ could also
be constructed via {\em Fourier-Mukai transform} in the following way.

Let $F$ be a sheaf on $\tproj$ and consider the diagram:
$$ \xymatrix{
& T \times \dproj \ar[dl]^{\Pi} \ar[dr]^{\hat{\Pi}} & \\
\tproj & & \dproj } $$
The Fourier-Mukai transform of $F$ is given by
$$ \Psi(F) = R\hat{\Pi}_* ( \Pi^*F \otimes {\cal P} ) $$
where ${\cal P}$ denotes the pullback of the Poincar\'e bundle
from $T\times\dual$ to $T \times \dproj$. If $F$ is torsion-free and
generically fibrewise semistable, then $\Psi(F)$ is a torsion sheaf on
$\dproj$. 

It is simple to show that if $F$ is locally-free and fibrewise 
semistable (as we have assumed throughout the paper), then $\Psi(F)$ is
supported exactly over the spectral curve $\overline{S}$, and the 
restriction to its support coincides with $\cal L$ \cite{JM}. 
Furthermore, it is also easy to see that 
$V=\pi_{1*}(R^1\hat{\Pi}_* ( \Pi^*F \otimes {\cal P}))$.

Therefore, the holomorphic version of the Nahm transform \cite{J1,J2} can be
seen as a Fourier-Mukai transform composed with Hitchin's correspondence. However, 
the Nahm transform (and the spectral construction of section \ref{instspec}) also 
contains some differential-geometric information (i.e. the instanton $A$, the
transformed connection $B$, and the spectral connection $\Gamma$) in addition 
to the holomorphic information encoded into the Fourier-Mukai transform.
 
Of course, such differential-geometric information is usually encoded into 
the holomorphic data in the form of a stability condition. Such condition is 
well-known for Higgs bundles \cite{H2}. For doubly-periodic instantons, the 
appropriate concept of stability for the corresponding instanton bundles is 
discussed in \cite{BJ}. It is less clear, though, what is the stability 
condition to be imposed on the spectral pairs $(\overline{S},\call)$; such 
question is addressed in \cite{JM}.


 \end{document}